\documentclass[reqno]{amsart}

\usepackage{amssymb}
\usepackage{amsmath}
\usepackage{mathrsfs}
\usepackage{graphicx}
\usepackage{fancybox}
\usepackage{tikz}

\newtheorem{theorem}{Theorem}

\newtheorem{proposition}[theorem]{Proposition}

\theoremstyle{definition}

\renewcommand{\geq}{\geqslant}

\newcommand{\Var}{{\textsc{Var}}}
\newcommand{\Form}{{\textsc{Form}}}

\newcommand{\eg}{\textit{e.g.}}
\newcommand{\ie}{\textit{i.e.}}

\title[Artificial precision and maximally consistent theories]{The  problem of artificial precision in theories of vagueness: a note on the r\^ole of maximal consistency}
\author{Vincenzo Marra}
\date{}
\address{Dipartimento di Matematica {\sl Federigo Enriques}, Universit\`a degli Studi di Milano, via Cesare Saldini 50, 20133 Milano, Italy}
\email{vincenzo.marra@unimi.it}

\begin{document}

\begin{abstract}
The problem of artificial precision is a major objection to any theory of vagueness based on real numbers as degrees of truth.  Suppose you are willing to admit that, under  sufficiently specified circumstances, a predication of ``is red''  receives a unique, exact number from the real unit interval $[0,1]$. You should then be committed to explain what is it that determines that value, settling for instance that my coat is red to
degree $0.322$ rather than $0.321$. In this note I revisit the problem in the important case of {\L}ukasiewicz infinite-valued propositional logic  that brings to the foreground the r\^ole of maximally consistent theories. I argue that the problem of artificial precision, as commonly conceived of in the literature, actually conflates two  distinct problems of a very different nature.
\end{abstract}

\maketitle
%
\section{}\label{s:intro}
 The monadic predicate $P(x):=$``$x$ is prime'', interpreted over the set of natural numbers $x\geq 1$, is (absolutely) \emph{precise}: its extension is  the set of prime numbers;
its anti-extension is the set of composite numbers;  each number either belongs to the extension of $P$ or to its anti-extension, but not to both; and in principle there is no issue as to whether a given number be prime or composite --- though in practice it may be impossible to ascertain which is the case for an astronomic instance of $x$. By contrast, the monadic predicate $R(x):=$``$x$ is red'', interpreted over the set of all objects, is (to some extent) \emph{vague}: its extension ought to be the set of all red objects; its anti-extension ought to be the set of all non-red objects; but it may not be  clear, even in principle, just which  objects do qualify as red, and which as non-red --- think of a peculiar tint at the borderline between red and pink.

Is there a logic of vague predicates --- or, for that matter, of vague propositions? Any definite answer, at present, is likely to be contentious. 
The philosophical literature on vagueness is relatively large; accounts of  the main competing theories may be found in \cite{Williamson1996, Keefe2000, shapiro, Smith2008}. One cluster of theories is based on the assumption that any instantiation of a vague predicate $R$ by a constant $c$ whose {\it denotatum} lies, intuitively, at the borderline between the extension of $R$ and its complement, is only true \emph{to a degree}. Thus if $c$ denotes my coat, and my coat is  of a peculiar tint at the borderline between red and pink, on  theories of this sort  the proposition $R(c):=$``My coat is red'' is to be considered neither true nor false, but rather true to some intermediate degree.  In this line of thought, a much stronger and yet popular assumption is that  the real unit interval $[0,1]$ embodies  ``degrees of truth''. See \cite[Chapter 4]{Williamson1996} for  an account of many-valued approaches to vagueness; \cite{Smith2008} is a recent proposal of a  theory of vagueness based on degrees of truth modelled by $[0,1]$.

\medskip Arguably, any such $[0,1]$-valued theory of vagueness faces the \emph{problem of artificial precision}.\footnote{The related but distinct problem of \emph{higher-order vagueness} (see \eg\ \cite[pp.\ 31--36]{Keefe2000}) will not be considered in this note.} Perhaps the first passage were the locution `artificial precision' was used in this connection is in
\cite[p.\ 443]{haack}:
\begin{quote}
[Fuzzy logic] imposes artificial precision [\ldots While] one is not obliged to
require that a predicate either definitely applies or definitely does not
apply, one is obliged to require that a predicate definitely applies to such-and-such, rather than to such-and-such other, degree (e.g.\ that a man 5ft\,10in tall belongs to \emph{tall} to degree 0.6 rather than 0.5).
\end{quote}
\noindent Tye \cite[p.\ 14]{tye} makes a similar point:

\begin{quote}One serious objection to [the many-valued approach] is that it really replaces vagueness with the most incredible and refined precision.
\end{quote}

\noindent In these terms, the objection is rather generic. It is not clear that a precise semantics for vague concepts is, \textit{per se}, a contradictory prospect. A sharper form of the objection, however, was put forth by Keefe  \cite[p.\ 47]{Keefe2000}, who identifies the source of the problem in our failure to see what could determine the correct value uniquely:

\begin{quote}
[T]he degree theorist's assignments impose precision in
a form that is just as unacceptable as a classical
true/false assignment. [\ldots] All predications of ``is red''
will receive a unique, exact value, but it seems
inappropriate to associate our vague predicate ``red''
with any particular exact function from objects to
degrees of truth. For a start, what could determine which
is the correct function, settling that my coat is red to
degree $0.322$ rather than $0.321$?
\end{quote}

\noindent  Smith \cite[p.\ 279]{Smith2008} endorses the objection in this form:

\begin{quote}
Intuitively, it is \emph{not} correct to say that there is one unique element of $[0,1]$ that 
correctly represents the degree of truth of `Bob is bald', with all other choices being incorrect. [\ldots]
 we have an affront to intuition [because] [w]e cannot see what could possibly \emph{determine} that the degree
of truth of `Bob is bald' is $0.61$ rather than $0.62$ or $0.6$ [...]
 \end{quote}

\smallskip  There is no question, I think, that the problem of artificial precision construed in this manner is a genuine objection to $[0,1]$-valued theories of vagueness; no such theory can get away without a plausible response to it. The purpose of this note is not to provide such a response. Rather more modestly, my aim is merely  to point out that the problem of artificial precision actually conflates    two  distinct issues of a very different nature: one falls within the realm of mathematical logic proper; the other belongs to enquiries into the semantics of vagueness. It seems to me that being clear about this distinction is an important preliminary to any treatment of the problem of artificial precision.

\section{}\label{s:basic}
The problem of artificial precision is only concerned with  propositions, because only propositions can receive a truth value. In order to focus on the essence of the objection, I will therefore restrict attention to propositional logic throughout the paper. However, even once it be agreed that truth values range in $[0,1]$, there remains much leeway to develop a formal system of many-valued propositional logic. The technical details of some of my arguments will depend on the specifics of the formal system under consideration; it is therefore important to be clear about the $[0,1]$-valued logic in question. I propose to concentrate on  \emph{{\L}ukasiewicz \textup{(}infinite-valued propositional\textup{)} logic}. This is a non-classical system going back to the 1920's, cf.\ the early survey  \cite[\S 3]{luktarski}, and its annotated  English translation  in \cite[pp.\ 38--59]{tarski}.
 The standard modern reference for  {\L}ukasiewicz logic is \cite{cdm}, while  \cite{mundicibis} deals
with topics at the frontier of current research. {\L}ukasiewicz logic can also be regarded as a member of a larger hierarchy of many-valued logics that was systematised by Petr H\'{a}jek in the late Nineties, cf.\ \cite{hajek}. Here I recall the basic notions.

\smallskip Let us fix once and for all the countably infinite set of propositional variables:
\[
\Var = \{X_1,X_2,\ldots,X_n,\ldots\}\,.
\]
Let us  write $\bot$ for the logical constant {\it falsum}, $\neg$ for the unary negation connective, and $\to$ for the binary implication connective. (Further derived connectives are introduced below.) The set $\Form$ of (well-formed) formul\ae\footnote{A set of conventions for omitting parentheses in formul\ae\ is usually adopted ($\bot$ is more binding than $\neg$, and $\neg$ is more binding than $\to$), and later extended to derived connectives. I do not spell the details here, as the conventions are analogous to the ones in classical logic, and are unlikely to cause confusion.} is defined exactly as in classical logic over the language $\{\bot,\neg,\to \}$. 

\smallskip The {\L}ukasiewicz calculus is defined by the five\footnote{In  \cite[Chapter 4]{cdm} the language has no logical constants, and consequently (A0) does not appear as an axiom. I prefer to explicitly have $\bot$ in the language, and thus I add {\it Ex falso quodlibet} to the standard axiomatisation.} axiom schemata
\begin{itemize}
\item[(A0)] $\bot \to \alpha$ \hfill (\textit{Ex falso quodlibet}.)
\item[(A1)] $\alpha\to(\beta\to\alpha)$ \hfill(\textit{A fortiori}.)
\item[(A2)] $(\alpha\to\beta)\to((\beta\to\gamma)\to(\alpha\to\gamma))$ \hfill(Implication is transitive.)
\item[(A3)] $((\alpha\to\beta)\to\beta)\to((\beta\to\alpha)\to\alpha)$ \hfill(Disjunction is commutative.)
\item[(A4)] $(\neg\alpha\to\neg\beta)\to(\beta\to\alpha)$ \hfill(Contraposition.)
\end{itemize}
with {\it modus ponens} as the only deduction rule. Provability is  defined exactly as in classical logic; $\vdash \alpha$ means that formula $\alpha$ is provable. I write $\mathscr{L}$ to denote {\L}ukasiewicz logic.

The logical constant \textit{verum} ($\top$), conjunction ($\wedge$), disjunction ($\vee$), and the biconditional ($\leftrightarrow$) are defined as in Table \ref{table:derivedconnectives}. From the definition of disjunction  one sees that (A3) indeed asserts the commutativity of disjunction. 
  Other common  derived connectives are reported in the same table, with their definition.
  \begin{table}[htf]
\begin{center}\begin{tabular}{|c|c|c|c|}
\hline {\bf Notation} & {\bf Definition} & {\bf Name}  & {\bf Idempotent} \\
\hline\hline
$\bot$& -- & {\it Falsum} &--\\
\hline
$\top$& $\neg \bot$  & {\it Verum}   &--\\
\hline
$\neg \alpha$ & -- & Negation  & --\\
\hline$\alpha\to \beta$ & -- & Implication  &--\\
\hline
$\alpha\vee \beta$ & $(\alpha\to \beta)\to \beta$ & (Lattice) Disjunction  & Yes\\
\hline
$\alpha\wedge \beta$ & $\neg(\neg\alpha \vee \neg\beta)$ & (Lattice) Conjunction  &Yes\\
\hline
$\alpha\leftrightarrow \beta$ & $(\alpha \to \beta)\wedge (\beta \to \alpha)$  &  Biconditional  &--\\
\hline $\alpha\oplus \beta$ & $\neg \alpha\to \beta$ &  Strong disjunction & No \\
\hline $\alpha\odot \beta$  & $\neg (\alpha\to \neg \beta)$ &  Strong conjunction & No\\
\hline $\alpha\ominus \beta$  & $\neg (\alpha\to \beta)$ &  But not, or Difference & -- \\
\hline \end{tabular} \smallskip 
\end{center}
\caption{Connectives in {\L}ukasiewicz logic.}\label{table:derivedconnectives}
\end{table}
  Some remarks are in order. Using the biconditional, one defines formul\ae\ $\alpha,\beta \in \Form$ to be \emph{logically equivalent} just in case $\vdash \alpha \leftrightarrow \beta$ holds. The connectives $\odot$ and $\oplus$ are then De Morgan dual: $\alpha \oplus \beta$ is logically equivalent to $\neg (\neg\alpha\odot\neg\beta)$, and $\alpha \odot \beta$ is logically equivalent to $\neg (\neg\alpha\oplus\neg\beta)$. These connectives, known as the \emph{strong disjunction} ($\oplus$) and \emph{strong conjunction} ($\odot$) of $\mathscr{L}$, play a central r\^ole both in H\'ajek's treatment of many-valued logics \cite{hajek}, and in Chang's algebraisation of $\mathscr{L}$ via MV-algebras \cite{cdm}. They are not idempotent,
 in the sense that  $\alpha \oplus \alpha$ and $\alpha$ are not logically equivalent: only the implication $\alpha\to\alpha\oplus\alpha$ is provable; dual considerations apply to $\odot$. Conjunction ($\wedge$) and disjunction ($\vee$) also are De Morgan dual, but they are  idempotent; in fact, they are sometimes called the \emph{lattice connectives} because they induce the structure of a distributive lattice in the algebraic semantics of $\mathscr{L}$. Finally, the connective $\ominus$ is the co-implication, \ie\ the dual to $\to$.

If $S\subseteq \Form$ is any set of formul\ae, one  writes $S\vdash \alpha$ to mean that $\alpha$ is provable  in {\L}ukasiewicz logic, under the additional set of  assumptions $S$. When this is the case, one says that $\alpha$ is a \emph{syntactic consequence} of $S$. Since each one of (A0--A4) is a principle of classical reasoning, and since \textit{modus ponens} is a classically valid rule of inference,  each formula provable  in $\mathscr{L}$ is a theorem of classical propositional logic. The converse is not true: most notably, the {\it tertium non datur} law, $\alpha\vee\neg \alpha$, is not provable in {\L}ukasiewicz logic; this  is one simple consequence of the completeness theorem to be recalled shortly. In fact, it can be shown that the addition of $\alpha\vee\neg \alpha$ as a sixth axiom schema to (A0--A4) yields classical logic.

\smallskip By a \emph{theory} in {\L}ukasiewicz logic one means any set of formul\ae\ that is closed under provability, \textit{i.e.}\ is deductively closed. For any $S \subseteq \Form$,  the smallest theory that extends $S$ exists: it is the \emph{deductive closure} $S^\vdash$ of $S$, defined by $\alpha \in S^\vdash$ if, and only if, $S\vdash \alpha$.  
A theory $\Theta$ is \emph{consistent} if $\Theta\not = \Form$, and \emph{inconsistent} otherwise; and it is \emph{maximal}, or \emph{maximally consistent}, if it is consistent, and whenever $\alpha\in \Form$ is such that $\alpha \not \in \Theta$, then $(\Theta\cup\{\alpha\})^\vdash=\Form$, \ie\ $\Theta \cup \{\alpha\}$ is inconsistent. 
A theory $\Theta$ is \emph{axiomatised by a set $S\subseteq \Form$ of formul\ae} if it so happens that $\Theta=S^\vdash$; and $\Theta$ is \emph{finitely axiomatisable} if $S$ can be chosen finite.

\smallskip
Let us now turn to the $[0,1]$-valued semantics. An \emph{atomic assignment},
or \emph{atomic evaluation}, is an arbitrary function $\overline{w}\colon \Var \to [0,1]$.  Such an atomic evaluation is uniquely extended to an \emph{evaluation} of all formul\ae, or \emph{possible world}, \textit{i.e.}\ to 
a function $w\colon \Form \to [0,1]$, via the compositional rules:
\begin{align*}
w(\bot)&=0\,,  \\
w(\alpha\to\beta)&=\min{\{1,1-(w(\alpha)-w(\beta))\}}\,,\\
w(\neg\alpha)&=1-w(\alpha)\,. 
\end{align*}
It follows by trivial computations that  the
formal semantics of derived connectives is the one reported in Table \ref{table:01connectives}.
\emph{Tautologies} are defined as those formul\ae\ that evaluate to $1$ under every evaluation.
\begin{table}[htf]
\begin{center}\begin{tabular}{|c|c|c|}
\hline {\bf Notation}    & {\bf Formal semantics} \\
\hline\hline
$\bot$& $w(\bot)=0$\\
\hline
$\top$  & $w(\top)=1$\\
\hline
$\neg \alpha$   & $w(\neg\alpha)=1-w(\alpha)$\\
\hline$\alpha\to \beta$   & $w(\alpha\to\beta)=\min{\{1,1-(w(\alpha)-w(\beta))\}}$\\
\hline
$\alpha\vee \beta$   & $w(\alpha\vee \beta)=\max{\{w(\alpha),w(\beta)\}}$\\
\hline
$\alpha\wedge \beta$   &$w(\alpha\wedge \beta)=\min{\{w(\alpha),w(\beta)\}}$\\
\hline
$\alpha\leftrightarrow \beta$ &$w(\alpha\leftrightarrow\beta)=1-|w(\alpha)-w(\beta)|$\\
\hline $\alpha\oplus \beta$  & $w(\alpha\oplus\beta)= \min{\{1,w(\alpha)+w(\beta)\}}$ \\
\hline $\alpha\odot \beta$   & $w(\alpha\odot\beta)= \max{\{0,w(\alpha)+w(\beta)-1\}}$\\
\hline $\alpha\ominus \beta$    & $w(\alpha\ominus\beta)= \max{\{0,w(\alpha)-w(\beta)\}}$ \\
\hline \end{tabular} \smallskip  
\label{table:01connectives}
\end{center}
\caption{ Formal semantics of connectives in {\L}ukasiewicz logic.}
\end{table}
Let us  write $\vDash \alpha$ to mean that the formula $\alpha\in \Form$ is a tautology. The relativisation of this concept to theories leads to the notion of semantic consequence. Let $S\subseteq \Form$ be any subset, and let $\Theta=S^\vdash$ be its associated theory. Given $\alpha\in \Form$, the assertion  $S\vDash \alpha$ states  that any evaluation $w\colon \Form\to [0,1]$ that satisfies $w(S)=\{1\}$ --- meaning that $w(\beta)=1$ for each $\beta \in S$ --- must also satisfy $w(\alpha)=1$. When this is the case, we say that $\alpha$ is a \emph{semantic consequence} of $S$. We write $S^{\vDash}$ for the set of semantic consequences of $S$.

\smallskip
It is an exercise to check that $\mathscr{L}$ enjoys the generalised validity theorem: for any $S \subseteq \Form$ and any $\alpha \in \Form$, if $S \vdash \alpha$ then $S \models \alpha$. On the other hand,   it is a  non-trivial theorem that $\mathscr{L}$ is complete\footnote{However, $\mathscr{L}$ fails strong completeness (\ie\ completeness for theories): there is a set $S\subseteq \Form$ and a formula $\alpha \in \Form$ such that $S\models\alpha$, but $S \not\vdash \alpha$; see \cite[4.6]{cdm}.} with respect to the many-valued semantics above: hence $\vdash \alpha$ if, and only if, $\models \alpha$, for any $\alpha \in \Form$. The first proof of this appeared in  \cite{rose_rosser}; see also  \cite[4.5.1 \& 4.5.2]{cdm}.

\smallskip
All of  the above can be adapted in the obvious manner to the finite set $\Var_n=\{X_1,\ldots, X_n\}$, in which case one speaks of {\L}ukasiewicz logic \emph{over $n$ \textup{(}propositional\textup{)} variables}, denoted $\mathscr{L}_n$. Although, strictly speaking, one should introduce fresh consequence relation symbols $\vdash_{n}$ and $\vDash_{n}$ for each  $\mathscr{L}_{n}$, I will avoid this pedantry and use  $\vdash$ and $\vDash$ in all cases; context will do the rest.  I will write $\Form_n$ for the set of   formul\ae\ whose propositional variables are contained in $\Var_n$. We will be mostly concerned with the apparently trivial case of $\mathscr{L}_{1}$.
\section{}\label{s:artificial}

\smallskip  What is the logical status of an assumption such as ``\,`VM is tall'\footnote{\label{fn:short}Not all vague predicates are alike. The predicate Tall\,($\cdot$) comes in pair with its opposite, Short\,($\cdot$), over the domain of all individuals, say; but Red\,($\cdot$) does not:  there is no colour term for  non-Red\,($\cdot$) in the visible spectrum. Although I will not argue the point in this note,  I believe that $\mathscr{L}$ cannot be an appropriate formal model of vague predicates such as  Red\,($\cdot$). Hence the shift from redness to tallness.} is true to degree $r\in [0,1]$\,''?
Let us begin with the remark of an eminent logician  \cite[p.\ 368]{hajekvague}.
\begin{quote}[Regarding the problem of artificial precision,]
[l]et us
comment that mathematical fuzzy logic concerns the possibility of sound inference,
surely not techniques of ascribing concrete truth degrees to concrete
propositions.
\end{quote}
\noindent The point H\'ajek is making here is by no means limited to mathematical fuzzy logic. Logic is exclusively concerned with the form of an argument, not with its content. Logic can teach us nothing (factual).

\smallskip Thus the assumption ``\,`VM is tall' is true to degree $r$\,'' has extra-logical content. In particular, it is a semantic assumption: it tells us that certain states of affairs (namely, those wherein ``VM is tall'' is true to a degree $\neq r$),
while perhaps logically consistent, are known (or assumed) not to be the case. In fact, it is reasonable to expect that the assumption ``\,`VM is tall' is true to degree $r$\,'' be \emph{maximally strong}, falling short only of the strongest, inconsistent assumption according to which everything is the case. For   the stronger an assumption is, the fewer  models it has, \ie\ the fewer are the possible worlds that are consistent with it. Now the assumption ``\,`\,VM is tall' is true to degree $r$\,''  leaves us with \emph{just one}\footnote{Assuming that all propositions under consideration are built from the single atomic one ``VM is tall'', and that the logic is truth-functional. Since, formally,   the discussion applies to $\mathscr{L}_{1}$, these assumptions are satisfied.} possible world consistent with it, namely, the one world in which VM is tall to degree \emph{exactly} $r$. This solitary possible world is the bare minimum we need to stay clear of the precipice of inconsistency.

\smallskip We can considerably sharpen these initial remarks. Let us focus on the  vague proposition
\begin{align}\tag{*}\label{t:p}
X_1:=\text{``VM is tall''}\,,\nonumber
\end{align}
formally modelled by the propositional variable $X_1$ in  $\mathscr{L}_{1}$.
There are two distinct situations.
\begin{enumerate}
\item[(S1)] All we know about $X_1$ is that it is a propositional variable. 
    \item[(S2)]  All we know about $X_1$ is that it is a propositional variable such that the possible worlds $w\colon \Form \to [0,1]$ that we are ready to admit in light of the intended interpretation (\ref{t:p}) of $X_1$ are precisely those satisfying $w(X_1)=r$, for a fixed real number $r\in [0,1]$.
\end{enumerate}
Failure to distinguish between (S1) and (S2) amounts to ignoring a set of extra-logical assumptions, namely, those encoded by (\ref{t:p}). Suppose first we are in situation (S1).  Consider any formula $\alpha(X_1) \in \Form_1$. Then, in our intended interpretation, $\alpha(X_1)$ is a statement about VM's tallness. But, given (S1), such a statement can be truthfully asserted if, and only if, it is a provable formula:  $\vdash \alpha(X_1)$. In other words, such statements coincide with analytic truths (relative to $\mathscr{L}_{1}${) which, by their very nature, \emph{are absolutely uninformative about  whoever's tallness}. Given  (S1) only, whatever statement one can truthfully assert about VM's tallness, one can equally truthfully assert about TW's thinness, and conversely. In situation (S1) the problem of artificial precision does not even make sense: there is no specific truth value to be puzzled about.\footnote{Compare  Haack's claim that ``[In fuzzy logic] one is obliged to require that a predicate definitely applies to such-and-such, rather than to such-and-such other, degree'' \cite[\textit{loc.\ cit.}]{haack}. It was just shown that the claim, if taken at face value, is  unwarranted.}

By contrast, it is  in situation (S2) that the problem of artificial precision arises. Now we are only concerned with evaluations
$w\colon \Form\to [0,1]$ that satisfy $w(X_1)=r$. Let us call such evaluations \emph{admissible}\footnote{There is just one such admissible evaluation in $\mathscr{L}_{1}$, of course. I am using the plural form in preparation for the forthcoming extension (S$_T$).} for our intended interpretation (\ref{t:p}) of $X_1$. Consider any formula $\alpha(X_1) \in \Form_1$. Then, on our intended interpretation, $\alpha(X_1)$ is a statement about VM's tallness. But, given (S2), it is no longer the case that $\alpha(X_1)$ can be truthfully asserted if, and only if, $\vdash \alpha(X_1)$. On the one hand, if $\vdash \alpha(X_1)$, then certainly  any admissible $w$ satisfies $w(\alpha(X_1))=1$ by the validity theorem, so that $\alpha(X_1)$ is indeed true on our intended interpretation (\ref{t:p}) --- but again, such analytic truths have nothing to do with my tallness. On the other hand, however, there will be formul\ae\ $\alpha(X_1)$ that are not provable in $\mathscr{L}_{1}$, but are such that $\alpha(X_1)$ \emph{can} be  truthfully asserted under (S2), precisely because  we restrict attention to admissible evaluations only. For a trivial example, assume $r=1$: then the formula $X_1$ can be truthfully asserted subject to (S2), simply because we restrict attention to the one possible world where $X_1$ indeed holds. In general, let $\Theta_r\subseteq \Form_1$ be the collection of all those formul\ae\ over the variable $X_1$ that may be truthfully asserted given (S2), that is, set
\begin{align}\tag{$\dag$}\label{t:Thetar}
\Theta_r=\{\alpha(X_{1})\in\Form_1 \ \mid \ w(\alpha(X_{1}))=1 \text{ whenever } w(X_1)=r\}\,,
\end{align}
where $w\colon \Form_{1}\to [0,1]$ ranges over all possible worlds.
Then those formul\ae\ in $\Theta_r$ that are not analytic truths \emph{are precisely the synthetic, factual truths about VM's tallness that the semantic assumption $w(X_1)=r$ entails, and that $\mathscr{L}_{1}$ is able to express syntactically}.
In other words, the set of formul\ae\ $\Theta_r$ attempts to encode the semantic assumption (S2) at the syntactic level, to within the formal linguistic resources afforded by  $\mathscr{L}_{1}$. 

\smallskip There are intermediate situations, of course.

\begin{enumerate}
\item[(S$_{T}$)]   All we know about $X_1$ is that it is a propositional variable such that the possible worlds $w\colon \Form \to [0,1]$ that we are ready to admit in light of the intended interpretation  (\ref{t:p})  of $X_1$ are precisely those satisfying $w(X_1)\in T$, for a fixed subset $T\subseteq [0,1]$.
\end{enumerate}
 If $T=[0,1]$, then (S$_{T}$) is (S1): we are imposing no restriction on possible worlds.  If $T=\{r\}$, then  (S$_{T}$) is (S2): the only possible world we are ready to admit is the one with $w(X_1)=r$.  We can adapt (\ref{t:Thetar}) to the intermediate situations in the obvious manner:
\begin{align}\tag{$\ddag$}\label{t:ThetaT}
\Theta_T=\{\alpha(X_{1})\in\Form_1 \ \mid \ w(\alpha(X_{1}))=1 \text{ whenever } w(X_1)\in T\}\,,
\end{align}
where $w\colon \Form_{1}\to [0,1]$ ranges over all possible worlds.

\smallskip
It turns out that  $\Theta_{T}$ as in (\ref{t:ThetaT}) is a theory in $\mathscr{L}_{1}$ for any choice of $T\subseteq [0,1]$:  the generalised validity theorem guarantees that the semantic assumption (S$_{T}$) is reflected into a deductively closed set of syntactic assumptions (\ref{t:ThetaT}). It is also clear by the very definition (\ref{t:ThetaT})  that  $\Theta_{T}$ is a superset of  $\Theta_{T'}$ whenever $T\subseteq T'$, for any two subsets $T,T'\subseteq [0,1]$. Thus
the (syntactic representation of the) assumption, say, ``\,`VM is tall' is true to degree $\geq \frac{1}{2}$\,'' is no stronger than the (syntactic representation of the) assumption ``\,`VM is tall' is true to degree $\frac{2}{3}$\,''. Once more, this suggests that a theory $\Theta_{r}$, for $r \in [0,1]$, should correspond to a maximally strong assumption.

\smallskip Let now $\mathscr{M}$ be the set of all
theories in $\mathscr{L}_{1}$ that can be written in the form (\ref{t:Thetar}). That is,
\begin{align}
\mathscr{M}=\{\Sigma\subseteq \Form_{1}\mid \text{ There exists } r \in [0,1] \text{ such that } \Sigma=\Theta_{r}\}\,.\nonumber
\end{align}
The next proposition\footnote{As mentioned,  $\Theta_r$ is deductively closed for any $r \in [0,1]$. Given $\alpha \in \Form_1$,
suppose $\Theta_{r} \vdash \alpha$. If $w_r\colon \Form_1 \to [0,1]$ is the
unique evaluation such that $w_r(X_1)=r$, then $w_r(\Theta_{r})=\{1\}$  by (\ref{t:Thetar}); since $\Theta_{r} \vdash \alpha$, then
 $w_r(\alpha)=1$ by the generalised
validity theorem for $\mathscr{L}_1$; hence $\alpha \in \Theta_r$, again by (\ref{t:Thetar}). Moreover, 
$\Theta_r$ is consistent: since $w_r(\bot)=0$ by the semantics of $\bot$, we have $\bot \not \in \Theta_r$ in light of (\ref{t:Thetar}).
It is harder to prove that $\Theta_r$ is maximally consistent, and that all maximally consistent theories are of this form. However, this is a standard result (essentially \cite[4.6.3 and 3.5.1]{cdm}).
} confirms the intuitions above about the members of $\mathscr{M}$.
\begin{proposition}\label{p:maxvalidity}  $\mathscr{M}$ is precisely the collection of all maximally consistent theories in $\mathscr{L}_1$.
\end{proposition}
\noindent 
This result  leaves open the possibility that different real numbers determine the same maximally consistent theory 
via (\ref{t:Thetar}), which leads to a key question.
\begin{enumerate}\item[(Q1)] Is the semantic assumption (S2)
precisely equivalent to the set of syntactic assumptions (\ref{t:Thetar})? More precisely, is the correspondence 
\[\tag{$\star$}\label{t:corr}
r \in [0,1] \ \longmapsto \ \Theta_{r} \in \mathscr{M}
\]
determined by definition
(\ref{t:Thetar}) a \emph{bijection} between the 
real unit interval $[0,1]$, and the
 set $\mathscr{M}$ of maximally consistent theories in $\mathscr{L}_{1}$?
\end{enumerate}
I note in passing that generalisations of this question make sense for arbitrary sets of propositional variables\footnote{In which case Proposition 1 and question (Q1) would be concerned with maximally consistent theories in $\mathscr{L}_{n}$ or $\mathscr{L}$.} (\ie\ for $\mathscr{L}_n$, $n > 0$ an integer, and for the whole of $\mathscr{L}$), and, in another direction, for arbitrary subsets $T\subseteq [0,1]$ of truth values and their associated  theories\footnote{In which case Proposition 1 and question (Q1) would be concerned with \emph{semisimple theories}, \ie\ those theories for which completeness holds; see \cite[4.6 and  3.6]{cdm}.} as in (\ref{t:ThetaT}). However, for the purposes of this note it will be enough to concentrate on (Q1). 

\smallskip Although the details will vary,\footnote{The main issue in generalizing (Q1) to other systems is that Proposition \ref{p:maxvalidity} most often fails, so that it is not enough to consider maximally consistent theories only. For example, in the important G\"odel-Dummett logic \cite[Chapter 4]{hajek}, $\mathscr{M}$ turns out to be exactly the collection of all prime theories in the one-variable fragment, where a theory $\Theta$ is \emph{prime} if it  proves either $\alpha \to \beta$ or $\beta\to \alpha$ for any two formul\ae\ $\alpha$ and $\beta$.} close analogues of (Q1) can be asked for virtually any $[0,1]$-valued logic $\mathscr{S}$, under rather weak assumptions. If the answer to (Q1) is negative then there is a discrepancy between our formal, $[0,1]$-valued semantics, and the expressive power of such a  logic $\mathscr{S}$. For if $\Theta_r=\Theta_s$ with $r\neq s \in [0,1]$, then  the logic $\mathscr{S}$ is not sufficiently expressive to discern between $r$ and $s$, so that it is legitimate to ask for further support to the claim that real numbers are the basis of a suitable formal semantics for $\mathscr{S}$. Indeed, proving that (Q1) has negative answer for a specific $[0,1]$-valued logic is a way of making precise the assertion that its $[0,1]$-valued semantics is  ``redundant''.  I do not wish to suggest that all redundant formal semantics (in the present sense) is useless. What I am implying, though, is that there is a sense in which a redundant semantics poses a challenge to  logicians:  \emph{some} justification for redundancy ought to be given, lest one incurs Occam's razor.\footnote{If $\mathscr{S}$ is such that (Q1) has negative answer, then the formal semantics of $\mathscr{S}$ violates  Leibniz's Identity of Indiscernibles: in deference to which, we ought to \emph{identify} real numbers (=possible worlds) $r,s \in [0,1]$ whenever they satisfy $\Theta_{r}=\Theta_{s}$, \emph{provided} all we are concerned with are those properties of $r$ and $s$ that are expressible within $\mathscr{S}$. If we insist not to
identify $r$ and $s$, then there must be distinguishing properties of these two possible worlds, not expressible with the linguistic resources of $\mathscr{S}$, that we nonetheless wish our formal semantics to record. A formal semantics strictly richer than the available language  is of course a  perfectly reasonable construct, but
it had better   result from an  explicit choice --- not from overlooking a negative answer to (Q1). I do not elaborate this line of thought further in  this note.} 
	 
\smallskip Such worries  need not concern us insofar as we are dealing with {\L}ukasiewicz logic:
\begin{proposition}\label{p:holder} The answer to \textup{(Q1)} is affirmative. That is, \textup{(}\ref{t:corr}\textup{)} is a bijection between $[0,1]$ and maximally consistent theories in $\mathscr{L}_{1}$.
\end{proposition}
\noindent  This  fact\footnote{Algebraically,  Proposition \ref{p:holder} amounts to the representation theorem for $1$-generated semisimple MV-algebras, see \cite[Chapter 3]{cdm}. Via Mundici's categorical equivalence between MV-algebras and lattice-ordered Abelian groups with a strong order unit \cite[7.1]{cdm}, this is equivalent to H\"older's theorem, for which the interested reader may consult \cite[\S 2.6]{bkw}.}  rests on Otto H\"older's theorem from 1901 that a totally ordered Ar\-chi\-me\-dean group embeds into the real numbers.
The mathematical details involved should not blind us to the remarkable conceptual content of Proposition \ref{p:holder}: \emph{the innocent-looking axioms \textup{(A0--A4)} of {\L}ukasiewicz logic characterise the real numbers}, in the sense that maximally consistent theories in $\mathscr{L}_{1}$ classify the elements of $[0,1]$. 

\section{} What does all this have to do with vagueness? Nothing:
(Q1) does not mention vagueness at all. {\it Mutatis mutandis}, (Q1)  is a question that can be asked about  virtually \emph{any}  logic with a $[0,1]$-valued semantics, and the answer will be independent of \emph{any} intended interpretation of that logic; briefly, (Q1) is a  question in mathematical logic.  

\smallskip But, having isolated (Q1) as worthy of attention, we can proceed to ask a second key question.
\begin{enumerate}
\item[(Q2)] Given $r \in [0,1]$, can we read off $\Theta_r$ as in (\ref{t:Thetar}) the synthetic truths about $X_1$ determined by (S2) in a form that is intelligible with respect to our intended interpretation (\ref{t:p})? 
\end{enumerate}
Vagueness ---  the intended interpretation  --- enters the picture at this stage: for  \textup{(Q2)} makes no sense in terms of the formal semantics alone.
We are given the truth value $r\in [0,1]$, and the corresponding  maximally consistent theory $\Theta_{r}$ induced in $\mathscr{L}_{1}$
by (\ref{t:Thetar}). We must provide, if we can, an interpretation of the formul\ae\ in $\Theta_{r}$ as assertions of vague propositions about VM's tallness,  that taken together amount to an explanation of what it means for VM to be tall to degree $r$. The answer to (Q2) is affirmative insofar as this interpretation and the ensuing explanation are convincing. In case $\Theta_r$ is finitely (hence singly) axiomatisable, and one such axiomatisation is known,\footnote{It can be proved that, given $r \in [0,1]$, $\Theta_{r}$ as in (\ref{t:Thetar}) is finitely axiomatisable in $\mathscr{L}_{1}$ if, and only if, $r$ is rational. Moreover, it is possible to exhibit an algorithm (for definiteness, a Turing machine) that, on input any rational number $r \in [0,1]$, outputs a formula $\alpha_{r}(X_{1})$ satisfying $\Theta_{r}=\{\alpha_{r}(X_{1})\}^{\vdash}$. Everything hinges on the theory of continued fractions and Schauder hats; see \cite[Chapter 3]{cdm}.}  we are thus concerned with a single formula $\alpha(X_1)$ such that $\Theta_r=\{\alpha(X_1)\}^\vdash$. We must provide a reading of $\alpha(X_1)$, as a single assertion about VM's tallness, that convincingly explains what it means to assume  ``\,`VM is tall' is true to degree $r$\,''. Observe that this explanation, whatever it is, need not (should not) mention $r$ itself --- nor, for that matter, any other number. This is  because $\alpha(X_1)$ is a formula in $\mathscr{L}_1$, and this logic has no truth constants other than \textit{verum} and \textit{falsum}.

\smallskip  I suggest here that (Q2) is a sharper formulation of the problem of artificial precision, as stated at the beginning of this note, for {\L}ukasiewicz logic --- and in fact, the specifics aside, for any aspiring $[0,1]$-valued logic of vagueness.

 \section{}
  So what is the answer to (Q2) for {\L}ukasiewicz logic?  New research\footnote{Since \cite{machina}, {\L}ukasiewicz logic has been widely discussed in the philosophical literature as a candidate for a logic of vagueness. More often than not, it has been rejected; cf.\ \eg\ \cite{Williamson1996, Keefe2000}. To the best of my knowledge, though, the quite specific question (Q2) has not been addressed.} is needed, I think, to say something defensible in this connection. Perhaps a glimpse of the difficulties involved may be caught if, by
 way of an epilogue to this note, 
 we  work our way up to the modestly non-classical truth value $\frac{1}{2}$.
 
 \smallskip The theory $\Theta_1$ is  axiomatised by the single formula $X_1$. In symbols,
 \[
 \Theta_1=\{X_1\}^\vdash\,.
 \]
In this case, the answer to (Q2) is that Your\footnote{Compare Bruno de Finetti's usage \cite{df1, df2} of the capitalised second person singular pronoun to stress that attributing a degree of probability to a (classical) proposition is a personal matter.} assumption ``\,`VM is tall' is true to degree $1$\,'' amounts to the fact that You are ready to assert that VM  is  clearly, indisputably tall.

\smallskip 
Similarly,  $\Theta_0$ is  axiomatised by the single formula $\neg X_1$:
 \[
 \Theta_0=\{\neg X_1\}^\vdash\,.
 \]
The answer to (Q2) is clear in this case, too. The assumption ``\,`VM is tall' is true to degree $0$\,'' amounts to the fact that You are ready to assert that VM  is  clearly, indisputably non-tall: in short,  short.\footnote{Cf.\ Footnote \ref{fn:short}.}

\smallskip What about $\Theta_{\frac{1}{2}}$? We have %
 \[
 \Theta_{\frac{1}{2}}=\left\{\,(\neg X_1\to X_1) \wedge (X_1\to \neg X_1)\,\right\}^\vdash\,.
 \]
Can we make sense of this? Quite generally (cf.\ Table \ref{table:01connectives}), to assert $\alpha \wedge \beta$ in $\mathscr{L}$ is  to assert $\alpha$ \emph{and} to assert $\beta$; and to assert $\alpha \to \beta$ is to assert that $\alpha$ is (clearly, indisputably) \emph{less true}\footnote{Caution: no circularity is involved in this passage. The objection of artificial precision can only be raised against theories that (i) have already committed to degrees of truth, and (ii) have embraced $[0,1]$,   or some other precisely specified structure, as a mathematical model for degrees of truth and their relationships.  The charge that we are here using a comparative notion of truth to explain artificial precision, without justifying the assumption  that truth  \emph{does} come in degrees, has therefore no force.
Similary, the problem of justifying why degrees of truth are modelled by the real numbers rather than, say, the octonions, may well be a problem --- there is no paucity of objections to (i--ii) in the literature --- but it is a different one.}
than $\beta$, or at the very most just as true. 

\smallskip So Your assumption that ``\,`VM is tall' is true to degree $\frac{1}{2}$\,'' amounts to the fact that You are ready to assert \emph{both} that VM  is  clearly, indisputably \emph{less} of a case of a short man, \emph{than} he is a case of a tall man, \emph{and} that VM  is  clearly, indisputably \emph{less} of a case of a tall man, \emph{than} he is a case of a short man.

\section*{Acknowledgements}
\noindent A preliminary version of parts of this paper was presented at the meeting {\it Epistemic Aspects of Many-Valued Logics}, held in Prague at the Institute of Philosophy of the Academy of Sciences of the Czech Republic, from the 13$^{\rm th}$ to the 16$^{\rm th}$ of September 2010. I am  grateful to the organisers, Timothy Childers, Christian Ferm\"{u}ller, and Ondrej Majer, for having given me a chance to present some of these ideas before an audience that included several philosophers who have thought deeply about vagueness. I am indebted to many participants for questions, discussions, and criticism that have been helpful in improving my initial ideas on the subject matter of this paper. In this connection, I should particularly like to thank Christian Ferm\"{u}ller, Colin Howson, Nicholas J.\ J.\ Smith,  and Timothy Williamson.

\bibliographystyle{amsalpha}
\providecommand{\bysame}{\leavevmode\hbox to3em{\hrulefill}\thinspace}
\providecommand{\MR}{\relax\ifhmode\unskip\space\fi MR }
\providecommand{\MRhref}[2]{%
  \href{http://www.ams.org/mathscinet-getitem?mr=#1}{#2}
}
\providecommand{\href}[2]{#2}

\end{document}